\documentclass{amsart}

\usepackage{amsmath}
\usepackage{mathtools}
\usepackage{amssymb}
\usepackage{enumerate}
\usepackage{slashed}
\usepackage{graphicx}
\usepackage{newlfont}
\usepackage{amsrefs}
\usepackage{comment}
\usepackage{mathtools}
\usepackage[top=25mm,bottom=25mm,right=25mm,left=25mm]{geometry}

\theoremstyle{plain}
\newtheorem{theorem}{Theorem}

\theoremstyle{definition}

\theoremstyle{remark}
\newtheorem{remark}{Remark}

\newcommand{\tr}{\operatorname{tr}}
\newcommand{\R}{\mathbb{R}}
\newcommand{\Z}{\mathbb{Z}}

\newcommand{\e}{\mathbf{e}}

\numberwithin{equation}{section}

\begin{document}

\title[Gravitational Solitons and Complete Ricci Flat Riemannian Manifolds]{Gravitational Solitons and Complete Ricci Flat Riemannian Manifolds of Infinite Topological Type}

\author[Khuri]{Marcus Khuri}
\address{Department of Mathematics\\
Stony Brook University\\
Stony Brook, NY 11794, USA}
\email{khuri@math.sunysb.edu}
\thanks{M. Khuri acknowledges the support of NSF Grant DMS-2104229, and Simons Foundation Fellowship 681443.}

\author[Reiris]{Martin Reiris}
\address{Centro de Matem\'{a}tica\\
Universidad de la Rep\'{u}blica\\
Montevideo, Uruguay}
\email{mreiris@cmat.edu.uy}

\author[Weinstein]{Gilbert Weinstein}
\address{Department of Mathematics \& Department of Physics\\
	Ariel University, Ariel, Israel 40700}
\email{gilbertw@ariel.ac.il}

\author[Yamada]{Sumio Yamada}
\address{Department of Mathematics\\
Gakushuin University\\
Tokyo 171-8588, Japan}
\email{yamada@math.gakushuin.ac.jp}
\thanks{S. Yamada acknowledges the support of JSPS KAKENHI Grant JP17H01091.}

\begin{abstract}
We present several new space-periodic solutions of the static vacuum Einstein equations in higher dimensions, both with and without black holes, having Kasner asymptotics. These latter solutions are referred to as gravitational solitons. Further partially compactified solutions are also obtained by taking appropriate quotients, and the topologies are computed explicitly in terms of connected sums of products of spheres.  In addition, it is shown that there is a correspondence, via Wick rotation, between the spacelike slices of the solitons and black hole solutions in one dimension less. As a corollary, the solitons give rise to complete Ricci flat Riemannian manifolds of infinite topological type and generic holonomy, in dimensions 4 and higher.
\end{abstract}
\maketitle

\begin{center}
\textit{This paper is dedicated to Demetrios Christodoulou. The work herein is inspired by the problem that he suggested to the third author, for his PhD thesis, almost four decades ago.}
\end{center}

\section{Introduction}\label{sec1}
\setcounter{equation}{0}
\setcounter{section}{1}

Configurations of multiple static and asymptotically flat vacuum black holes typically exhibit conical singularities and possibly other pathologies. These cone angles are necessary to keep the gravitational equilibrium and are therefore usually interpreted as repulsive forces. A detailed analysis of this phenomenon in the axisymmetric setting was originally carried out by Bach and Weyl \cite{bachweyl1} for the 4-dimensional case, proving in particular the existence of an angle defect along any piece of the axis joining two black holes. In higher dimensions, a thorough investigation was initiated by Emparan and Reall \cite{ER}. We know that in general, regular configurations of two or more black holes are ruled out in 4 dimensions by Bunting and Masood-ul-Alam's theorem \cite{BuntMasood1}. This result also extends to higher dimensions \cites{GibbonsShiromizu1,HwangS1}, where the Schwarzschild-Tangherlini solution realizes the only asymptotically flat static black hole in a vacuum. Without the assumption of asymptotic flatness, however, this type of rigidity breaks down even in dimension 4. Indeed, partly motivated by the balancing of multiple static electro-vacuum black holes in the Majumdar-Papapetrou spacetime, Myers \cite{myers19871} constructed regular 4-dimensional static vacuum spacetimes in which an infinite number of Schwarzschild black holes are aligned in a periodic fashion along an axis of symmetry. These same solutions were later rediscovered by Korotkin and Nicolai in \cite{korotkinnicolai1}. The Myers-Korotkin-Nicolai solutions are asymptotically Kasner, and play an integral role in an extended version of static black hole uniqueness given by Peraza and Reiris \cite{PerazaReiris1}.  In \cite{myers19871}, it was conjectured that these space-periodic vacuum configurations can be generalized to higher dimensions, perhaps with black holes of nontrivial topology. We confirmed this to be the case in \cite{KhuriWeinsteinYamadaJHEP1}, by producing a variety of examples of 5-dimensional space-periodic static vacuum solutions with combinations of the sphere $S^3$ and ring $S^1\times S^2$ horizon cross-sectional topologies.

The methods of \cite{KhuriWeinsteinYamadaJHEP1} also led to the discovery of 5-dimensional vacuum solitons. Here, a \textit{gravitational soliton} refers to a nontrivial, globally static and geodesically complete spacetime. Although such horizonless soliton solutions are well-known features of supergravity theories \cite{Kundurietal1}, where nontrivial topology is supported by magnetic flux supplied through Maxwell fields, an asymptotically flat stationary vacuum spacetime which is geodesically complete must be Minkowski space. This latter statement is a classical result of Lichnerowicz \cite{Lichnerowicz1} in dimension 4, and in general can be established as a corollary of the rigidity portion of the positive mass theorem \cites{schoenyau19791,SchoenYau20171,witten19811}, together with Stokes' theorem and the Komar expression for mass. This no-soliton result essentially holds even without the assumption of asymptotic flatness in 4 dimensions, in the sense that solitons of this dimension are always covered by the Minkowski spacetime \cite{anderson1}*{Theorem 0.1}. In dimensions greater than four, vacuum solitons must have constant lapse \cites{Chen,Reiris}, and therefore such spacetimes factor into a pure product of time with a complete Ricci flat Riemannian manifold (a Cauchy hypersurface).   

The examples constructed in \cite{KhuriWeinsteinYamadaJHEP1}*{Theorem 2} have a time slice topology homeomorphic to an infinite connected sum $\#^{\infty} \> S^2\times S^2$, and therefore have infinite second Betti number. These solutions of the static vacuum equations admit Kasner asymptotics, are bi-axisymmetric, geodesically complete,  and space-periodic. By taking quotients, we also obtain solitons on $M_k^4 \setminus (B^2 \times T^2 )$ where $B^2$ and $T^2$ are the 2-ball and 2-torus respectively, and $M_k^4$ is either $S^4$ or $\#^{k} \> S^2\times S^2$ depending on whether $k=0$ or $1\leq k<\infty$.

The purpose of the present work is as follows. We will extend the results of \cite{KhuriWeinsteinYamadaJHEP1} to obtain space-periodic vacuum solitons in all dimensions greater than 3, and will then confirm that the static potential for these solutions must in fact be constant, thus the time slices yield complete Ricci flat Riemannian manifolds. These manifolds are simply connected but of infinite topological type, admit Kasner asymptotics and generic holonomy, and are of cohomogeneity-two via a torus action. Moreover, these manifolds are periodic in one direction, allowing for further solutions by taking discrete quotients; we are able to compute the topology of such quotients in terms of connected sums of products of spheres. In addition, it will be shown that there is a correspondence, induced through Wick rotation, between the space-periodic Ricci flat manifolds arising from the solitons, and space-periodic black hole solutions in one dimension less. Previous examples of complete Ricci flat manifolds of infinite topological type have been constructed in 4-dimensions by Anderson-Kronheimer-LeBrun \cite{AKL} using the Gibbons-Hawking ansatz, and by Goto \cite{Goto} in dimensions $4m$ with a $T^m$ symmetry using the hyper-K\"{a}hler quotient method; these were later studied further by Hattori \cite{Hattori}. The approach of Goto was also generalized by Dancer-Swann \cite{DancerSwann} to produce hypertoric manifolds of infinite topological type, which were analyzed in more detail in \cite{Dancer}. All of these previous works fit within the hyper-K\"{a}hler context, and in particular have dimensional restrictions as well as special holonomy. By contrast, the Ricci flat manifolds that we produce are derived from an entirely different source, and consequently they exhibit distinct properties. More precisely, the construction arises naturally as a byproduct from our study of the stationary vacuum multi-axisymmetric Einstein equations in higher dimensions \cites{KWY0,KWY1,KWY2,KWY3,KhuriWeinsteinYamadaJHEP1,KWetal}. These equations reduce to a study of singular harmonic maps from $\mathbb{R}^3$ into a nonpositively curved symmetric space, although, in the static case with some additional restrictions, solutions may be found by the superposition of Green's functions with concentration along intervals of the $z$-axis. The arrangement of such intervals and the choice of Green's functions determine the so called `rod structure' of the solution, which may be used to prescribe the topology and other aspects of the spacetime. Typically conical singularities are produced in this process, however, we show that by choosing certain periodic arrangements of the rods these singularities are relieved.

\section{Background and Setup}\label{sec2}
\setcounter{equation}{0}
\setcounter{section}{2}

Let $(\mathbf{M}^{n+3},\mathbf{g})$, $n\geq 1$ be the domain of outer communication of a stationary $n$-axisymmetric $(n+3)$-dimensional spacetime, that is, it admits $\mathbb{R}\times U(1)^n$ as a subgroup of its isometry group. Under reasonable hypotheses \cite{hollands2011}, the orbit space $\mathbf{M}^{n+3}/[\mathbb{R}\times U(1)^n]$ is homeomorphic to the right half plane $\{(\rho,z)\mid \rho>0\}$. In this setting, the vacuum Einstein equations reduce to an axisymmetric harmonic map, with domain $\mathbb{R}^3\setminus \{z-\text{axis}\}$ parameterized by the cylindrical coordinates $(\rho,z,\varphi)$, and target symmetric space \mbox{$SL(n+1,\mathbb{R})/SO(n+1)$}. The $z$-axis boundary of the orbit space is decomposed into an exhaustive sequence of intervals called \textit{rods}, each of which is defined by a particular isotropy subgroup of $U(1)^n$. We will label these intervals by $\{\Gamma_l\}_{l\in I}$ for some index set $I$ which may be infinite, and will divide the set of rods into two types, namely \textit{axis rods} and \textit{horizons rods}. Each axis rod $\Gamma_l$ is defined by the vanishing of a linear combination $v^i_l \partial_{\phi^i}$ of the generators $\partial_{\phi^i}$, $i=1,\ldots,n$ of the $U(1)^n$ symmetry, where the nonzero vector $\mathbf{v}_l=(v_l^1,\ldots,v_l^n)\in\mathbb{Z}^n$ consists of relatively prime integers so that $\mathrm{gcd}\{v_l^1,\ldots,v_l^n\}=1$, and is called the \textit{rod structure} of $\Gamma_l$; the coordinates $\phi^i$ on $T^n$ have period $2\pi$. Thus, each axis rod structure defines a 1-dimensional isotropy subgroup $ \mathbb{R}/\mathbb{Z} \cdot \mathbf{v}_l \subset \mathbb{R}^n / \mathbb{Z}^n \cong T^n$ for the action of $T^n$ on points that lie over $\Gamma_l$.
On the other hand, a horizon rod $\Gamma_h$ is an interval of the $z$-axis where no closed-orbit Killing field degenerates, that is $\mathbf{v}_h=0$, but where $|\partial_t+\Omega^i\partial_{\phi^i}|$ vanishes with $\partial_t$ denoting the stationary Killing field, and with $\Omega_i$, $i=1,\ldots, n$ representing the horizon angular velocities. A point in the orbit space at which two neighboring axis rods intersect is referred to as a \textit{corner}, and since two (linearly independent) rotational Killing fields vanish there, the total space over this point is an $(n-2)$-torus; the intersection point of an axis rod with a horizon rod is called a \textit{pole}. In order to avoid orbifold singularities at corner points, neighboring axis rod structures $\mathbf{v}$, $\mathbf{w}$ are required to satisfy \cite{KWetal}*{Section 3} the \textit{admissibility condition} $\mathrm{Det}_2(\mathbf{v},\mathbf{w})=\pm 1$, where the second determinant divisor is defined by
\begin{equation}
\mathrm{Det}_2(\mathbf{v},\mathbf{w})=\mathrm{gcd}\{Q_{\mathbf{j}}\mid \mathbf{j}=(j_1,j_2)\in\mathbb{Z}^2, 1\leq j_1<j_2\leq n\}
\end{equation}
with the determinant $Q_{\mathbf{j}}$ arising from the $2\times 2$ minor constructed from the matrix having columns $\mathbf{v}$, $\mathbf{w}$ by choosing the $j_1$, $j_2$ rows. The collection of rods and associated rod structures completely determines the topology of horizons and the domain of outer communication, see \cites{KWY2,KWetal}. For instance, if the two axis rods bordering a horizon rod have rod structures $\mathbf{e}_i$, $\mathbf{e}_j$, where $\mathbf{e}_i$ is the standard basis element of $\mathbb{Z}^n$ which has a 1 in the $i$th slot and zeros elsewhere, then the horizon topology is $S^3 \times T^{n-2}$ for $i\neq j$ and $S^1 \times S^2 \times T^{n-2}$ for $i=j$.

The stationary $n$-axisymmetric vacuum Einstein equations \cite{maison} reduce to solving the following harmonic map system
\begin{align}\label{23jjikjgh}
\begin{split}
&\Delta f_{ij}-f^{ab}\nabla^{k}f_{ia}\nabla_{k}f_{jb}
+f^{-1}\nabla^{k}\omega_{i}\nabla_{k}\omega_{j}=0,\\
&\Delta\omega_{i}-f^{ab}\nabla^{k}f_{ia}\nabla_{k}\omega_{b}
-f^{ab}\nabla^{k}f_{ab}\nabla_{k}\omega_{i}=0,
\end{split}
\end{align}
where $\Delta$ is the $\mathbb{R}^3$-Laplacian, $F=(f_{ij})$ is an $n\times n$ symmetric matrix which is positive definite away from the axes, $f=\det F$, and $\omega=(\omega_1,\ldots,\omega_n)^{t}$ is the set of twist potentials associated with the $U(1)^n$ symmetry. These quantities parameterize the symmetric space target through a $(n+1)\times (n+1)$ symmetric, positive definite, unimodular matrix
\begin{equation}
\Phi=
\begin{pmatrix}
f^{-1}&  -f^{-1}\omega_i\\
-f^{-1}\omega_i  & f_{ij}+f^{-1}\omega_i\omega_j
\end{pmatrix},
\quad
i,j=1,..,n.
\end{equation}
Furthermore, the spacetime metric on $\mathbf{M}^{n+3}$ can be constructed from these quantities and expressed in Weyl-Papapetrou coordinates by
\begin{equation}\label{aofig}
\mathbf{g}=-f^{-1}\rho^2 dt^2+e^{2\alpha}(d\rho^2+dz^2)
+f_{ij}(d\phi^{i}+\beta^{i}dt)(d\phi^{j}+\beta^{j}dt).
\end{equation}
Note that this shows an interpretation of rod structures as vectors $\mathbf{v}_l$ lying in the 1-dimensional kernel of the matrix $F$ at an axis rod $\Gamma_l$. The functions $\alpha$ and $\beta^i$ may be obtained by quadrature \cite{shiromizu}, more precisely
they can be found by integrating the equations
\begin{equation}\label{fs00}
\beta^i_{\rho}=\rho f^{-1}f^{ij}\omega_{j,z},\quad\quad\quad
\beta^i_{z}=-\rho f^{-1}f^{ij}\omega_{j,\rho},
\end{equation}
and
\begin{align}\label{a}
\begin{split}
\alpha_\rho =&
\frac\rho8 \left[ (\log f)_\rho^2 - (\log f)_z^2 + \tr F^{-1}F_\rho F^{-1}F_\rho
- \tr F^{-1}F_z F^{-1}F_z
- \frac4\rho (\log f)_\rho + \frac2{f}F^{-1}(\omega_\rho^2-\omega_z^2) \right], \\
\alpha_z =& \frac\rho4 \left[ (\log f)_\rho(\log f)_z + \tr F^{-1}F_\rho F^{-1}F_z -\frac2\rho (\log f)_z + \frac2{f} F^{-1}\omega_{\rho}\omega_z \right],
\end{split}
\end{align}
where we have used the notation $F^{-1}\omega_{\rho}\omega_z:=\omega_{\rho}^t F^{-1}\omega_z$. The integrability conditions for \eqref{fs00} and \eqref{a} correspond to the harmonic map equations \eqref{23jjikjgh}. If $\Gamma$ denotes the union of all axis rods, then the relevant harmonic map $\Phi:\mathbb{R}^3\setminus\Gamma\rightarrow SL(n+1,\mathbb{R})/SO(n+1)$ is singular along the axes, and its asymptotics encode the rod structures and values of the potentials at these points. Boundary conditions (prescribed asymptotics) are therefore imposed on the axes in order to achieve the desired rod structures, and the potentials $\omega$ are assigned to be constants $\mathbf{c}_l\in\R^n$ on each axis rod $\Gamma_l$, in such a manner to guarantee that the values of the constants agree on consecutive axis rods. Hence, the potential constants can only change across horizon rods, and the difference determines the horizon angular momenta. The relevant existence theory is studied in \cite{KWetal}.

A solution to the singular harmonic map problem gives rise to an $n$-axisymmetric stationary vacuum spacetime, with prescribed rod structures and horizon angular momenta. However, it is possible that conical singularities form on the axes when assembling the spacetime metric \eqref{aofig} from the harmonic map. The conical singularity at an interior point $(0,z_0)$ along the axis rod $\Gamma_{l}$, having rod structure $\mathbf{v}_l$, is quantified by the angle defect $\theta\in(-\infty,2\pi)$ arising from the 2-dimensional cone formed by the orbits of $v^j\partial_{\phi^j}$ over the line $z=z_0$ in the orbit space. This value may be computed from the expression 
\begin{equation}\label{ffeehjhr}
\frac{2\pi}{2\pi-\theta}=\lim_{\rho\rightarrow 0}\frac{2\pi\cdot\mathrm{Radius}}{\mathrm{Circumference}}
=\lim_{\rho\rightarrow 0}
\frac{\int_{0}^{\rho}\sqrt{e^{2\alpha}}}{\sqrt{f_{ij}v^{i}v^{j}}}
=\lim_{\rho\rightarrow 0}
\sqrt{\frac{\rho^2 e^{2\alpha}}{f_{ij}v^{i}v^{j}}}.
\end{equation}
A conical singularity is absent if the angle defect vanishes.
It is routine to check that with the aid of a change from polar to Cartesian coordinates, this condition is necessary and sufficient for the smooth extendibility of the metric across the axis, assuming that analytic regularity has been established. Moreover, analytic regularity allows for a well-defined notion of \textit{logarithmic angle defect} $b_l=\log\left(\frac{2\pi}{2\pi-\theta}\right)$ associated with the axis rod $\Gamma_l$, since the angle defect must then be constant on each axis rod \cite{harmark}. The conical singularity on $\Gamma_l$ is referred to as an \textit{angle deficit} if $b_l>0$, and an \textit{angle surplus} if $b_l<0$.

\section{Statement of Results}\label{sec3}
\setcounter{equation}{0}
\setcounter{section}{3}

We will now restrict attention to static $n$-axisymmetric vacuum spacetimes. This requires the vanishing of twist potentials $\omega_i=0$, $i=1,\ldots,n$, and significantly simplifies the harmonic map equations \eqref{23jjikjgh}. Note, however, that the equations are still nonlinear. In order to make contact with a linear system, we may impose a further ansatz that restricts the metric along the torus fibers to be given as a diagonal matrix function
\begin{equation}
F=\mathrm{diag}(e^{u_1},\ldots,e^{u_n}).
\end{equation}
Observe that with these assumption, the vacuum Einstein equations reduce to finding $n$ harmonic functions $u_i$ on $\mathbb{R}^3 \setminus\Gamma$, and the spacetime metric takes the form
\begin{equation}\label{aofig1}
\mathbf{g}=-\rho^2 e^{-\sum_{i=1}^{n}u_i}dt^2+e^{2\alpha}(d\rho^2+dz^2)+\sum_{i=1}^{n}e^{u_i} \left(d\phi^i \right)^2.
\end{equation}
In this setting the axes can only exhibit the rod structures $\mathbf{e}_i$, $i=1,\ldots,n$ from the standard basis of $\mathbb{Z}^n$. For an axis rod $\Gamma_l$ having the rod structure $\mathbf{e}_l$, we find that the corresponding logarithmic angle defect is given by
\begin{equation}\label{b1-b2}
b_l =  \lim_{\rho \rightarrow 0}  \left(\log \rho +\alpha-\frac{1}{2}u_l\right)\quad\text{ on }\quad\Gamma_l.
\end{equation}
In what follows, the functions $u_i$ will be constructed as a sum of Green's functions such that they are periodic in the $z$-direction. The function $\alpha$, which is obtained by quadrature from \eqref{a}, will then be shown to also possess the same periodicity, yielding the desired space-periodic static vacuum spacetimes. More precisely, we will say that such solutions are \textit{space-periodic} if the group $\mathbb{Z}$ acts by isometries through translations in the $z$-direction of the Weyl-Papapetrou coordinate system. The spacetimes that we discuss are \textit{asymptotically Kasner}, meaning that the metric asymptotes to the Kasner form
\begin{equation}\label{aktrho}
\mathbf{g}\thicksim -q_0 dt^2 +q_1 d\tau^2 +\tau^{2p_0} dz^2 +\sum_{i=1}^{n} \tau^{2p_i} \left(d\phi^i \right)^2,
\end{equation}
where $q_0 ,q_1>0$ are constants and the exponents satisfy the Kasner conditions $\sum_{i=0}^n p_i =\sum_{i=0}^n p_i^2 =1$. Geometric regularity of the solutions is established by eliminating the possibility of conical singularities along the axes. This will be achieved by utilizing the degrees of freedom arising from addition of constants to the $u_i$ and $\alpha$.

\begin{theorem}\label{thm1}
For each $n\geq2$, there is a regular $n$-axisymmetric static vacuum soliton spacetime $(\mathbf{M}^{n+3},\mathbf{g})$ which is space-periodic and asymptotically Kasner. The rod structure is periodic with fundamental period $\e_1,\dots,\e_n $. Furthermore, these spacetimes are simply connected and of infinite topological type, in that the codimension-three Betti number is infinite, $b_{n}(\mathbf{M}^{n+3})=\infty$.
\end{theorem}

The staticity of the spacetimes implies that topologically $\mathbf{M}^{n+3}=\mathbb{R}\times M^{n+2}$, where $M^{n+2}$ represents a constant time slice. The periodicity allows for the taking of quotients by subgroups of $\mathbb{Z}$, to obtain further solutions with partially compactified time slice topology $\tilde{M}^{n+2}$. Since the rod structure period of the quotients contains the fundamental basis for $\mathbb{Z}$, these spaces will also be simply connected. Moreover, $\tilde{M}^{n+2}$ still admits an effective $T^n$-action, and therefore we may apply the classification results of Oh and Orlik-Raymond \cites{oh,oh1,orlikraymond} to compute the topology of these slices. Up to spatial dimension 6, within the spin category, closed simply connected manifolds admitting a cohomogeneity-two torus action are connected sums of products of spheres. By choosing appropriate periodic configurations of rod structures, beyond those treated in Theorem \ref{thm1}, we are able to produce solitons whose topology involves each type of spherical product appearing in the classification list. Moreover, for dimensions greater than or equal to 6, the topology associated with the basic rod structure sequence is also known to be a connected sum of products of spheres, by the work of McGavran \cite{MG}*{Theorem 3.4} and additionally \cite{Bosio}*{Theorem 6.3}, \cite{BP}*{Theorem 4.6.12}. In the statements below, note that cases $(iii)$ and $(iv)$ have a partial overlap.

\begin{theorem} \label{existence}
There exist regular $n$-axisymmetric static vacuum soliton spacetimes \mbox{$(\mathbb{R}\times \tilde{M}^{n+2},\mathbf{g})$} which are asymptotically Kasner, and such that $\tilde{M}^{n+2}$ admits the following topologies with associated rod structures.
\begin{enumerate}[(i)]
\item For $n=2$, the spatial slice is homeomorphic to either $S^4 \setminus (B^2 \times T^2)$ or $S^2 \times S^2 \setminus (B^2 \times T^2)$,
with rod structure periods $\mathbf{e}_1, \mathbf{e}_2$ or $\mathbf{e}_1, \mathbf{e}_2, \mathbf{e}_1, \mathbf{e}_2$ respectively.
    
\item For $n=3$, the spatial slice is homeomorphic to either $S^5 \setminus (B^2 \times T^3)$ or $S^2 \times S^3 \setminus (B^2 \times T^3)$,
with rod structure periods $\mathbf{e}_1, \mathbf{e}_2, \mathbf{e}_3$ or $\mathbf{e}_1, \mathbf{e}_2, \mathbf{e}_1, \mathbf{e}_3$ respectively.

\item  For $n=4$, the spatial slice is homeomorphic to either $S^3 \times S^3 \setminus (B^2 \times T^4)$ or $\left[(S^2 \times S^4) \# 2(S^3 \times S^3)\right] \setminus (B^2 \times T^4)$, with rod structure periods $\mathbf{e}_1, \mathbf{e}_2, \mathbf{e}_3, \mathbf{e}_4$ or $\mathbf{e}_1, \mathbf{e}_2, \mathbf{e}_1, \mathbf{e}_3 ,\mathbf{e}_4$ respectively.

\item For $n\geq 4$, the spatial slice is homeomorphic to 
\begin{equation}\label{khgjqohgh}
\left[\overset{n-3}{\underset{k=1}{\#}}k {n-2 \choose k+1}  S^{2+k}\times S^{n-k}\right]\setminus (B^2 \times T^n), 
\end{equation}
with rod structure period $\mathbf{e}_1,\ldots,\mathbf{e}_n$. 
\end{enumerate}
\end{theorem}

Soliton solutions are devoid of horizons. However, the solitons exhibited in Theorems \ref{thm1} and \ref{existence} can be used to produce static vacuum black hole solutions in one dimension less. In particular, by choosing an angular coordinate $\phi^i$, we may Wick rotate the Riemannian spatial slice so that $\partial_{\phi^i}$ becomes timelike. Thus, the time slices of the solitons become black hole spacetimes in which the $\mathbf{e}_i$-axis rods transform into horizon rods. It can then be shown that the resulting spacetimes are vacuum. The horizon topologies resulting from this process are restricted to the product of a sphere or ring, with tori. We note that the original space-periodic black hole solution found by Myers-Korotkin-Nicolai \cites{korotkinnicolai1,myers19871} arises from this process, applied to the soliton of Theorem \ref{thm1} with $n=2$.

\begin{theorem}\label{wick}
Let $(M^{n+2},g)$, $n\geq 2$ be a time slice of any soliton produced in Theorems \ref{thm1} or \ref{existence}, and pick an angular coordinate $\phi^i$, $1\leq i\leq n$. Then the Wick rotation $\phi^i \rightarrow \sqrt{-1}t$ applied to the Weyl-Papapetrou form of the metric transforms this manifold into a regular black hole solution $(\mathbf{M}^{n+2},\mathbf{g})$ of the $(n-1)$-axisymmetric static vacuum Einstein equations, with Kasner asymptotics. 
\end{theorem}

The converse of this result also holds, in that given a space-periodic $(n-1)$-axisymmetric static vacuum solution, Wick rotation of the time coordinate into an angular variable produces a periodic $n$-axisymmetric Ricci flat Riemannian manifold of the same dimension, which can then serve as the time slice of a soliton in one dimension higher. The key observation needed to establish this statement, as well as to prove Theorem \ref{wick}, is that the static potential for the solitons produced in Theorems \ref{thm1} and \ref{existence} must be constant. A consequence of this fact is that the time slices of these solitons are complete Ricci flat Riemannian manifolds, and thus yield new examples of Riemannian Einstein metrics. We will say that such Riemannian manifolds are asymptotically Kasner if the metric asymptotes to the time slice of the Kasner metric in \eqref{aktrho}.

\begin{theorem} \label{lapse}
The time slice $(M^{n+2},g)$ of any soliton from Theorems \ref{thm1} or \ref{existence} is a complete Ricci flat Riemannian manifold, admitting a cohomogeneity-two torus action, and with Kasner asymptotics. Furthermore, those arising from Theorem \ref{thm1} are simply connected and of infinite topological type, in that the codimension-two Betti number is infinite, $b_{n}(M^{n+2})=\infty$.
\end{theorem}

Previous examples of complete Ricci flat Riemannian manifolds of infinite topological type have been found within the hyper-K\"{a}hler context, as discussed in the introduction. By contrast, those arising from the solitons above are of generic holonomy. This is proven, with the help of the Ambrose-Singer Theorem, by analyzing the structure of the curvature tensor in the asymptotic end.

\begin{theorem} \label{holonomy}
The complete Ricci flat Riemannian manifolds of Theorem \ref{lapse} are of generic holonomy.
\end{theorem}


\section{Soliton Existence}\label{sec4}
\setcounter{equation}{0}
\setcounter{section}{4}

In this section we establish existence of the solitons asserted in Theorems~\ref{thm1} and \ref{existence}. The arguments follow closely those of~\cite{KhuriWeinsteinYamadaJHEP1}. For simplicity, we will assume that in the rod configurations all rod lengths are equal. In a remark at the end of the section, we will describe how this requirement may be relaxed with some restrictions.

\subsection{Theorem \ref{thm1} existence}
\label{4.1}

Let $L>0$, and divide the $z$-axis of $\R^3$ into segments $\Gamma_l$, $l\in\Z$ such that each is of length $L/n$. The initial step consists of constructing $n$ axially symmetric harmonic functions $u_i$, $i=1,\dots n$ on $\R^3$, which are $L$-periodic in $z$ and asymptote to $2\log\rho$ near the rods $\Gamma_{nl+i}$, $l\in\Z$. For completeness we present this construction here. Consider the Green's function for a uniform charge distribution along an interval $I=[a,b]$ within the $z$-axis, namely 
\begin{equation}
U_I=\log(r_a-z_a)-\log(r_b-z_b),
\end{equation}
where 
\begin{equation}
r_a = \sqrt{\rho^2+(z-a)^2}, \quad\quad\quad z_a = z-a. 
\end{equation}
Note that this function satisfies the following properties
\begin{equation}\label{fjgh}
U_I<0,\quad\quad\quad U_I \thicksim 2\log\rho\quad\text{ near }\quad I,\quad\quad\quad
U_I=(a-b)/r + O(r^{-2})\quad\text{ as }\quad r\to\infty.
\end{equation}
We may then form the potentials
\begin{equation}\label{UIP}
u_i = \lim_{m\to\infty} \left( \sum_{l=-m}^m U_{\Gamma_{nl+i}} + \frac{2}{n} \log m \right),\quad\quad i=1,\ldots,n.
\end{equation}
For any $(\rho,z)$ with $\rho>0$, and each large $l$, we have by \eqref{fjgh} that $U_{\Gamma_{nl+i}}(\rho,z)\thicksim \frac{1}{nl}$. Thus, the additional term $\frac{2}{n}\log m$ renormalizes the divergent series of harmonic functions to produce a finite harmonic limit away from the axis. Near each $\Gamma_{nl+i}$, $l\in\Z$ the asymptotics for these functions is $u_i\sim 2\log\rho$, so that according to \eqref{aofig1} the rod structure of $\Gamma_{nl+i}$ is $\mathbf{e}_i$. The harmonic functions $u_i$ are $L$-periodic by construction, and it can be shown as in \cite{KhuriWeinsteinYamadaJHEP1}*{Example 1} by expanding in cylindrical harmonics (modified Bessel functions) that 
\begin{equation}\label{jiqughh}
u_i =\frac{2}{n} \log\rho +c_i +O(e^{-a_i \rho}) \quad\text{ as }\quad \rho\rightarrow\infty,\quad\quad i=1,\ldots,n,
\end{equation}
for some constants $a_i >0$ and $c_i$. 

With the harmonic functions $u_i$ in hand, we may form the spacetime metric \eqref{aofig1}.
In order to aid with the arguments showing that the function $\alpha$ is $L$-periodic, it should be observed that the sum of these functions reduces to an explicit analytic expression, namely
\begin{equation} \label{constant-lapse}
\sum_{i=1}^n u_i = 2\log\rho +c,
\end{equation}
where $c$ is the sum of $c_i$. To see this, observe that $\sum_{i=1}^n u_i -2\log\rho -c$ is a harmonic function that is uniformly bounded, $L$-periodic, and tends to zero as $\rho\rightarrow\infty$. One can then apply a version of the maximum principle \cite{Weinstein}*{Lemma 8} on horizontal strips in the $\rho z$-plane to conclude that this function must vanish identically. In fact, it can be shown that $c=-2\log (2L)$ with the following alternative justification of \eqref{constant-lapse}. By choosing the origin of the $z$-axis appropriately we have
\begin{align}
\begin{split}
\sum_{i=1}^{n} u_i =&\lim_{m\rightarrow\infty}\left(U_{\left[-\tfrac{(2m+1)}{2}L,\tfrac{(2m+1)}{2}L\right]}
+2\log m\right)\\
=&\lim_{m\rightarrow\infty}\left(U_{\left[-\tfrac{(2m+1)}{2}L,\tfrac{(2m+1)}{2}L\right]}
+2\log \left((2m+1)L\right)\right)-2\log (2L)\\
=& 2\log\rho-2\log (2L).
\end{split}
\end{align}

We will now consider the periodicity of $\alpha$. Note that with \eqref{constant-lapse}, the quadrature equations \eqref{a} simplify to
\begin{equation} \label{alpha}
\alpha_\rho = \frac{\rho}{8} \, \sum_{i=1}^n \left(u_{i,\rho}^2-u_{i,z}^2 \right) - \frac{1}{2\rho}, \quad\quad\quad
\alpha_z = \frac{\rho}{4} \,\sum_{i=1}^n u_{i,\rho} u_{i,z}.
\end{equation}
The second of these equations may be used to show that $\alpha$ is $L$-periodic, that is, 
if we fix $\rho$ then
\begin{equation}\label{period}
\int_P \alpha_z\, dz = 0
\end{equation}
where the integration is carried out over one period. Indeed, since all the $u_i$'s are periodic, we can carry out the integration for each term $u_{i,\rho} u_{i,z}$ on a different period. Clearly, if we take the period centered at the midpoint of $\Gamma_i$, then $u_i$ is even with respect to this midpoint, which implies that $u_{i,\rho}u_{i,z}$ is odd and~\eqref{period} follows.

\begin{figure}
	\centering
	\includegraphics[width=10cm]{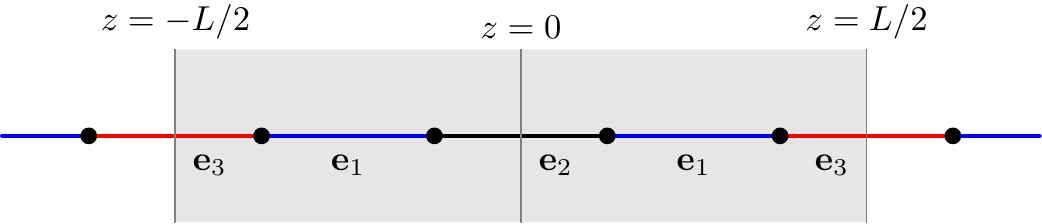}
	\caption{Periodicity of $\alpha$ in the $\e_1,\e_2,\e_1,\e_3$ soliton.}
	\label{fig-6-d}
\end{figure}

Consider now the issue of conical singularities. Recall that the logarithmic angle defect on axis rod $\Gamma_i$ is given by
\begin{equation}
b_i =\lim_{\rho\rightarrow 0} \left(\log\rho +\alpha -\frac{1}{2}u_i\right), \quad\quad i=1,\ldots,n.
\end{equation}
Since the $u_i$ are only determined up to the addition of constants, we are free to adjust these  constants to ensure that $b_i =0$, $i=1,\ldots,n$. By periodicity, all axes are then free of conical singularities. Furthermore, note that in a neighborhood of a point in the interior of $\Gamma_i$, we have that $\bar{u}_i=u_i -2\log\rho$ is smooth. This follows from the fact that $\bar{u}_i$ is uniformly bounded and harmonic in this neighborhood away from the axis by construction, and \cite{Weinstein}*{Lemma 8} may be used to show that it can be extended smoothly across the axis. Regularity of the spacetime metric \eqref{aofig1} is then established with arguments analogous to those of \cite{KhuriWeinsteinYamadaJHEP1}*{Section 5.1}; regularity at the corners is treated similarly.  This completes the proof of existence for Theorem \ref{thm1}.

\subsection{Theorem \ref{existence} existence} 

All of the solitons described in this theorem arise from space-periodic solitons by taking a quotient. In particular, if $(M^{n+2},g)$ is the time slice of a space-periodic soliton with  fundamental period $L$, consider the discrete isometry group $\mathbb{Z}$ acting on $M^{n+2}$ by $z \mapsto z+ L$. This action is clearly properly discontinuous, and hence the quotient $\tilde{M}^{n+2}$ is a Riemannian manifold. Conical singularities will be absent from $\tilde{M}^{n+2}$ if they are absent from the cover $M^{n+2}$. Thus, it suffices to show that regular space-periodic solitons exist with the given rod configurations.

Case $(i)$ of this result is given in \cite{KhuriWeinsteinYamadaJHEP1}*{Theorem 2}. Moreover, the first set of rod structures within case $(ii)$ is covered by the previous theorem, and so we now consider the second set of rod structures $\mathbf{e}_1 , \mathbf{e}_2, \mathbf{e}_1, \mathbf{e}_3$. As in Section \ref{4.1}, the $z$-axis may be divided into rod intervals of equal length having an $L$-periodic configuration, with a fundamental domain consisting of this sequence of rod structures. Furthermore, $L$-periodic harmonic functions $u_1$, $u_2$, $u_3$ can be constructed that respect the given rod structure configuration. It remains to confirm the periodicity of $\alpha$, and to balance any conical singularities on the axes. To establish periodicity of $\alpha$, note that by translating, we can set the midpoint of an axis rod corresponding to $\e_2$ to be $z=0$. Then all three potentials $u_1$, $u_2$, and $u_3$ are even with respect to the line $z=0$, see Figure~\ref{fig-6-d}. It then follows as before that each term $u_{i,\rho} u_{i,z}$ is odd in the expression for $\alpha_z$, and therefore \eqref{period} holds, showing that $\alpha$ is periodic. In order to relieve any conical singularities, we must arrange for the logarithmic angle defects to vanish on the four axis rods of a fundamental domain. Note that there are three free constants arising from the potentials $u_i$, and thus we can immediately balance three of the axis rods in a fundamental domain, say those associated with the $\mathbf{e}_2$, $\mathbf{e}_3$ rods and one of the $\mathbf{e}_1$ rods. Furthermore, the rod structure configuration in Figure~\ref{fig-6-d} clearly admits an involutive symmetry defined by reflection across the line $z=0$, and this is manifest in the functions $u_i$, $i=1,2,3$ and $\alpha$. It follows that the logarithmic angle defect must also vanish for the remaining $\mathbf{e}_1$ rod in the fundamental domain, since it coincides with the image of the balanced $\mathbf{e}_1$ rod under the involution. We then have regular $L$-periodic solitons devoid of conical singularities in case $(ii)$.

\begin{figure}
	\centering
	\includegraphics[width=10cm]{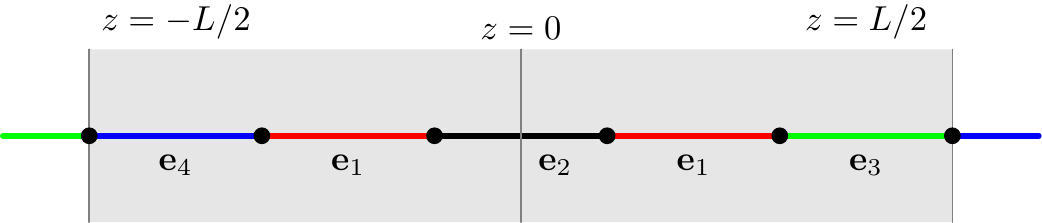}
	\caption{Periodicity of $\alpha$ in the $\e_1,\e_2,\e_1,\e_3,\e_4$ soliton.}
	\label{fig-7-d}
\end{figure}

The first set of rod structures within case $(iii)$ is covered by Theorem \ref{thm1}, and so we now consider the second set of rod structures $\mathbf{e}_1 , \mathbf{e}_2, \mathbf{e}_1, \mathbf{e}_3, \mathbf{e}_4$. As before, the $z$-axis is divided into rod intervals of equal length with this sequence of rod structures defining an $L$-periodic configuration, and four $L$-periodic harmonic functions $u_1,\ldots,u_4$ are constructed that respect the rod structures. To prove periodicity of $\alpha$, we again choose the midpoint of an axis rod corresponding to $\e_2$ to be $z=0$, as in Figure~\ref{fig-7-d}. It follows that $u_3(\rho,-z)=u_4(\rho,z)$ (up to addition of constants) and $u_1$, $u_2$ are even functions in $z$. As a consequence $\alpha_z$ is odd since
\begin{align}
\begin{split}
u_{3,\rho}(\rho,-z) u_{3,z}(\rho,-z)=&-u_{4,\rho}(\rho,z)u_{4,\rho}(\rho,z),\\
u_{i,\rho}(\rho,-z) u_{i,z}(\rho,-z)=&-u_{i,\rho}(\rho,z)u_{i,\rho}(\rho,z),\quad\quad i=1,2.
\end{split}
\end{align}
Thus \eqref{period} holds, confirming that $\alpha$ is periodic.
In order to relieve any conical singularities, observe that there are four free constants arising from the potentials $u_i$, which can be used to balance four of the axis rods in a fundamental domain, say those associated with the $\mathbf{e}_2$, $\mathbf{e}_3$, $\mathbf{e}_4$ rods, and one of the $\mathbf{e}_1$ rods. Moreover, the rod structure arrangement in Figure~\ref{fig-7-d} clearly admits an involutive symmetry across the line $z=0$, and this is transmitted to the functions $u_1$ and $\alpha$ as described above. 
In analogy with the previous example, it follows that the logarithmic angle defect must also vanish for the remaining $\mathbf{e}_1$ rod in the fundamental domain, due to this symmetry. We then have regular $L$-periodic solitons devoid of conical singularities in case $(iii)$. 

\begin{figure}
	\centering
	\includegraphics[width=8.5cm]{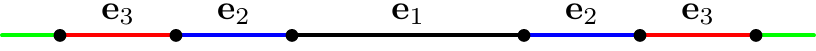}
	\caption{Symmetry of the fundamental domain with non-equal axis rod lengths.}
	\label{non-equal}
\end{figure}

\begin{remark} \label{extension}
We note that the existence arguments presented above can be extended to non-equal rod lengths and other configurations with some limitations. As an example we describe here an extension with one rod having a different length than the rest in a fundamental domain, see Figure~\ref{non-equal}. To see that the same proof carries over, notice that if the origin of the $z$-axis is placed at the center of the $\mathbf{e}_1$ rod, then the potential functions $u_i$ and $\alpha$ are even functions with respect to $z$. As before, this implies that $\alpha$ is periodic. Furthermore, we may use the free constants associated with the $u_i$ to balance the conical singularities on the axis rods within the fundamental domain which intersect the positive $z$-axis. By the involutive symmetry, the remaining axis rods in the fundamental domain will also be balanced, yielding a regular soliton. This generalizes the procedure of case $(ii)$. Similarly, the procedure of case $(iii)$ also admits a generlization to the non-equal rod lengths regime. In this situation, the involutive symmetry is more complicated in that it arises not only from a reflection in the domain space, but also involves composition with a map that interchanges two harmonic functions $u_i$, $u_j$ with $i\neq j$. These latter maps may be viewed as isometries of the target space when discussing the harmonic map formulation of this problem. Again, the symmetry allows us to balance the second half of the axis rods in a fundamental domain, after balancing the first half using the free constants associated with the potentials. 
\end{remark}

\section{Soliton Topology and Asymptotics}\label{sec5}
\setcounter{equation}{0}
\setcounter{section}{5}

In this section we will establish the topological claims of Theorems \ref{thm1} and \ref{existence}, as well as the Kasner asymptotics. First observe that the solitons of these theorems are all simply connected, since the integer span of their rod structures is  $\mathbb{Z}^n$ \cite{KWetal}*{Theorem 7.1}. Below we will verify the specific topologies of the quotient solitons, and show that the space-periodic solitons are of infinite topological type.

\subsection{Topology in Theorem \ref{existence}}

Within the space-periodic solitons, the fundamental domain of rod structures is of length $L$, which identifies a strip $-L/2\leq z\leq L/2$ in the orbit space. After identification to obtain the partially compactified solitons, this strip becomes a punctured disc $B^2 \setminus \{0\}$ with boundary circle that is divided into the various rod structures making up the fundamental domain. In order to analyze the topology, we may fill-in the asymptotic end (represented by the puncture) with $B^2 \times T^n$ to obtain a compact manifold $\hat{M}^{n+2}$, whose orbit space is the whole disc with boundary circle dividing into rod structures. It then suffices to describe the topology of the compactified manifold $\hat{M}^{n+2}$.

\begin{figure}
\centering
\includegraphics[width=3cm]{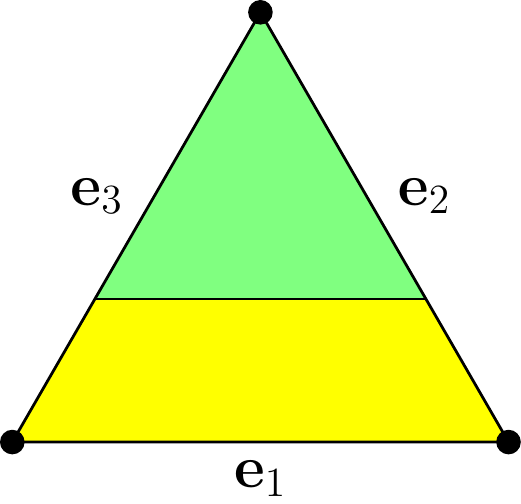}
\caption{Time slice topology of the $\e_1,\e_2,\e_3$ soliton.}
\label{triangle}
\end{figure}

\begin{figure}
\centering
\includegraphics[width=3cm]{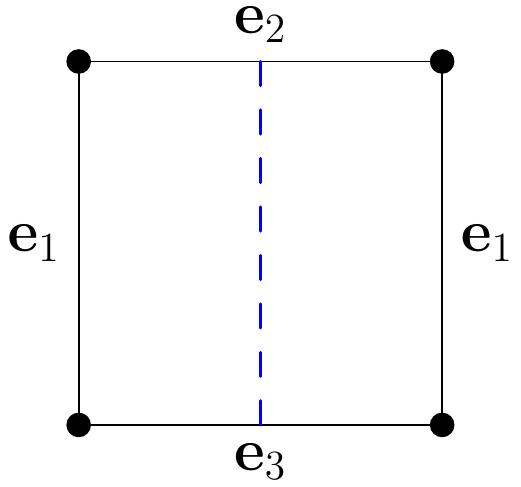}
\caption{Time slice topology of the $\e_1,\e_2,\e_1,\e_3$ soliton.}
\label{square1}
\end{figure}

Case $(i)$ is treated in \cite{KhuriWeinsteinYamadaJHEP1}*{Theorem 2}, so consider case $(ii)$. In the first example the orbit space consists of a disc with boundary rod structures
$\e_1,\e_2,\e_3$, see Figure~\ref{triangle}. The  green shaded region lifts to  $B^4\times S^1$, a solid ring, and the yellow shaded region lifts to $S^3\times B^2$. The two are glued along an $S^3\times S^1$, hence yielding an $S^5$. Next, consider the example with rod structures  $\e_1,\e_2,\e_1,\e_3$, see Figure~\ref{square1}. Similarly to the previous case, each vertical line, such as the dashed blue line, lifts to an $S^3\times S^1$. The $S^1$ collapses to a point at the left and right edges, hence yielding an $S^3\times S^2$.

\begin{figure}
\centering
\includegraphics[width=3cm]{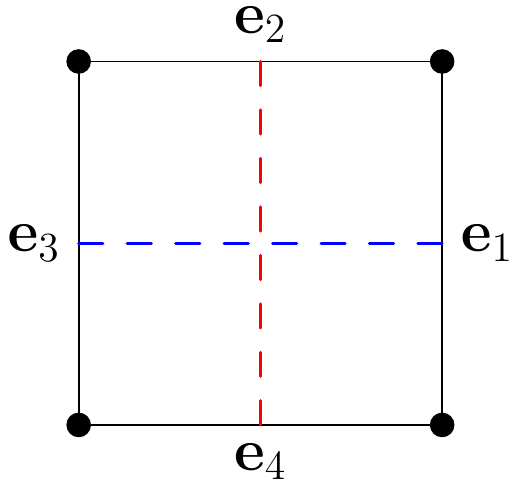}	
\caption{Time slice topology of the $\e_1,\e_2,\e_3,\e_4$ soliton.}
\label{square}
\end{figure}

We now examine case $(iii)$ starting with the rod structure sequence $\e_1,\e_2,\e_3,\e_4$, as illustrated in Figure~\ref{square}. Each of the vertical lines, such as the vertical dashed red line in the figure, lifts to an $S^3\times T^2$. Note that each of the horizontal lines, such as the dashed blue line, also lifts to an $S^3\times T^2$. Thus the topology corresponding to this rod diagram is clearly a Cartesian product. The $T^2$, in the $S^3\times T^2$ vertical slice degenerates along a different generator at the left and right edges, hence yielding an $S^3\times S^3$. 

\begin{figure}
\centering
\includegraphics[width=3cm]{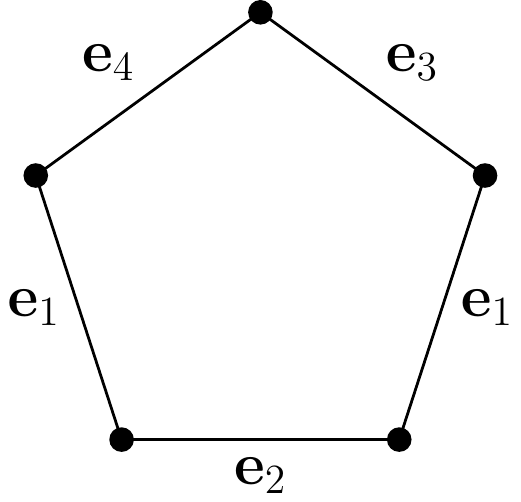}
\caption{Time slice topology of the $\e_1,\e_2,\e_1,\e_3,\e_4$ soliton.}
\label{pentagon}
\end{figure}

Next, we classify the topology of case $(iii)$ having rod structures $\e_1, \e_2, \e_1, \e_3, \e_4$, see Figure~\ref{pentagon}. This example is different in that it is not apparent how to obtain the topology via the methods used above. Instead, we will appeal to classification results of \cites{oh,oh1}. For convenience when comparing with these references, rename the rod structures by
\begin{equation}
\mathbf{e}_1\rightarrow \mathbf{e}_3, \quad\quad \mathbf{e}_2 \rightarrow\mathbf{e}_4,\quad\quad\mathbf{e}_3 \rightarrow \mathbf{e}_1,\quad\quad\mathbf{e}_4 \rightarrow\mathbf{e}_2,
\end{equation}
then the circle boundary of the orbit space has rod structure sequence
$\e_3, \e_1, \e_2, \e_3, \e_4$. Observe that the circle action on $\hat{M}^6$ associated with the subgroup of $T^4$ generated by the $2\mathbf{e}_1 +\mathbf{e}_2 -\mathbf{e}_4$ generator, is free. The 5-dimensional quotient manifold $\hat{M}^{6}/\!\!\sim$ then admits an effective $T^3$ action, and has a disc orbit space with rod structures $\e_3, \e_1, \e_2, \e_3, 2\mathbf{e}_1 +\mathbf{e}_2$. Furthermore, the proof of \cite{oh1}*{Theorem 5.5} shows that this manifold is spin, and because it is also simply connected and has five rods we find that $\hat{M}^{6}/\!\!\sim \text{ }\!\!\!\!\cong\# 2 \left(S^3 \times S^2 \right)$. By employing the Whitney product formula, and using that $\hat{M}^{6}/\!\!\sim$ is spin, it follows that $\hat{M}^6$ is spin. Since this manifold is simply connected and has five rods, the classification of \cite{oh}*{Theorem 1.1} implies that $\hat{M}^6 \cong (S^2 \times S^4) \# 2(S^3 \times S^3)$.

Lastly, consider case $(iv)$ in which $n\geq 4$, and the orbit space for $\hat{M}^{n+2}$ is determined by the basic sequence of rod structures $\e_1,\ldots,\e_n$. In \cite{MG}*{Theorem 3.4}, it is shown that a closed simply connected $(n+2)$-dimensional manifold having an effective $T^n$ action, with exactly $n$ rods on the orbit space boundary, must be the connected sum of products of spheres in \eqref{khgjqohgh}. By changing coordinates in the torus fibers (see \cite{KWetal}*{Lemma 3.3}), such a sequence of rod structures may be transformed into the basic sequence. Thus, we conclude that $\hat{M}^{n+2}$ has the topology arising from \eqref{khgjqohgh} by filling in the asymptotic end. We note that although McGavran's paper \cite{MG} contains an error in Theorem 3.6, as pointed out by Oh \cite{oh1}, Theorem 3.4 of \cite{MG} is not affected and has been generalized in \cite{Bosio}*{Theorem 6.3} (see also \cite{BP}*{Theorem 4.6.12}).

\subsection{Infinite topological type}
\label{5.2}

Here we show that the codimension-two Betti number $b_n(M^{n+2})$ is infinite, for the time slices of solitons produced in Theorem \ref{thm1}. Since we are only concerned with topology in this subsection, a new metric may be introduced on $M^{n+2}$ which is not necessarily free of conical singularities, namely
\begin{equation}
g_0=d\rho^{2}+dz^{2}+\sum_{i=1}^{n}e^{u_{i}}\left(d\phi^{i}\right)^{2}
\end{equation}
where $u_i$, $i=1,\ldots,n$ are as in Section \ref{4.1}. Let $\Gamma_{nl+i}$ be a rod on the $z$-axis with associated Green's function $U_{\Gamma_{nl+i}}$. This function is harmonic on $\mathbb{R}^3 \setminus \Gamma_{nl+i}$ with respect to the flat metric. However, we may also view $U_{\Gamma_{nl+i}}$ as an $n$-axisymmetric function on $M^{n+2}$ which is harmonic with respect to $g_0$, in light of \eqref{constant-lapse}. Therefore the $n$-form on $M^{n+2}$ given by \begin{equation}
\iota_{\eta_i}\star dU_{\Gamma_{nl+i}}   
=e^{c/2}\rho \left(\partial_{\rho} U_{\Gamma_{nl+i}} dz\wedge d\phi^{2}\wedge \cdots\hat{d\phi^i}\cdots\wedge d\phi^{n} - \partial_{z} U_{\Gamma_{nl+i}} d\rho\wedge d\phi^{2}\wedge\cdots\hat{d\phi^i}\cdots\wedge d\phi^{n} \right)
\end{equation} 
is closed by Cartan's formula, where $\star$ is the $g_0$-Hodge star and $\iota_{\eta_i}$ denotes interior product with $\eta_i=\partial_{\phi^i}$. Moreover, the asymptotics at the axis \eqref{fjgh} show that this form is smooth on $M^{n+2}$. Next, note that the axis rod $\Gamma_{nl+i}$ 
lifts to an embedded $n$-cycle $\Sigma_{nl+i}\cong S^{3}\times T^{n-3}$ in $M^{n+2}$ for $n\geq 3$, with $\Sigma_{nl+i}\cong S^2$ when $n=2$. Furthermore
\begin{equation}
\int_{\Sigma_{nl+i}}\iota_{\eta_i}\star dU_{\Gamma_{nl+i}}=2(2\pi)^{n-1}e^{c/2}|\Gamma_{nl+i}|\neq 0,\quad\quad\quad
\int_{\Sigma_{j}}\iota_{\eta_i}\star dU_{\Gamma_{nl+i}}=0,\quad\quad j\neq nl+i,
\end{equation}
where $|\Gamma_{nl+i}|$ is the length of the rod. It follows that each $[\Sigma_{nl+i}]$ represents a distinct generator of the homology group $H_{n}(M^{n+2};\mathbb{Z})$, yielding the desired conclusion.

\subsection{Asymptotics}

We will now confirm the asymptotics of the solitons constructed in Theorems \ref{thm1} and \ref{existence}. Recall that the Kasner metric on $\mathbb{R}^{n+1,1}$ takes the form
\begin{equation}
g_{K}=-dt^2 +\sum_{i=0}^n t^{2p_i} \left(dx^i\right)^2,
\end{equation}
and that this metric satisfies the vacuum Einstein equations exactly when the Kasner conditions hold:
\begin{equation}\label{alkjfhpag}
\sum_{i=0}^n p_i =1,\quad\quad\quad\quad \sum_{i=0}^n p_i^2=1.
\end{equation}
The solitons that we produce have metrics of the form~\eqref{aofig1}, where the asymptotics of the coefficients as $\rho\rightarrow\infty$ are determined by
\begin{equation}
u_i\thicksim A_i\log\rho, \quad
\quad\quad\alpha\thicksim C\log\rho,
\end{equation}
with $A_i>0$ and 
\begin{equation}\label{faihgiojoqiohgh}
C=\frac{1}{8}\sum_{i=1}^n A_i^2 - \frac12.
\end{equation}
It follows that
\begin{equation}
\mathbf{g}\thicksim -4L^2 dt^2+\rho^{2C}(d\rho^2 +dz^2) + \sum_{i=1}^{n}\rho^{A_i} \left(d\phi^i\right)^2,
\end{equation}
where we have used \eqref{constant-lapse}.
Since $\mathbf{g}$ solves the vacuum Einstein equations,
the powers of $\rho$ in the above expression satisfy the Kasner conditions. To see this more explicitly, set $\tau=\rho^{C+1}$ and observe that $C+1>0$ as well as
\begin{equation}
\mathbf{g}\thicksim - 4L^2 dt^2+\frac{1}{(C+1)^2} d\tau^2+\tau^{\frac{2C}{C+1}}dz^2 
+ \sum_{i=1}^n\tau^{\frac{A_i}{C+1}}\left(d\phi^i\right)^2.
\end{equation}
We can then verify, using $\sum_{i=1}^{n} A_i = 2$ from \eqref{jiqughh}, that the Kasner conditions~\eqref{alkjfhpag} hold for any values of $A_i$, so long as $C$ is given by~\eqref{faihgiojoqiohgh}. More precisely
\begin{equation}
\frac{C}{C+1} +\sum_{i=1}^{n} \frac{A_i}{2(C+1)} = 1,\quad\quad
\left(\frac{C}{C+1}\right)^2 + \sum_{i=1}^{n} \frac{A_i^2}{4(C+1)^2} = \frac{C^2+2C+1}{(C+1)^2} = 1.
\end{equation}
Therefore, the solutions of Theorems \ref{thm1} and \ref{existence} are asymptotically Kasner, with the role of `time' being played by the spatial variable $\rho$ when the metric is considered in the Kasner context.

\section{Wick Rotation, Riemannian Einstein Metrics, and Holonomy}\label{sec6}
\setcounter{equation}{0}
\setcounter{section}{6}

In this section we will establish Theorems \ref{wick}, \ref{lapse}, and \ref{holonomy}. Let $(M^{n+2},g)$, $n\geq 2$ denote the time slice of a space-periodic soliton arising from Theorem \ref{thm1} or \ref{existence}; recall that the solutions of Theorem \ref{existence} derive from quotients of  space-periodic solitons.  Choose an angular coordinate $\phi^i$ on the torus, and perform the Wick rotation $\phi^i \rightarrow\sqrt{-1} t$ 
which transforms the spacelike part of the metric~\eqref{aofig1} into
\begin{align}
\begin{split}
\mathbf{g} =& -e^{u_i} dt^2 + e^{2\alpha}(d\rho^2+dz^2)+ \sum_{j\ne i} e^{u_j} \left(d\phi^j\right)^2 \\
=&-\rho^2 e^{-\sum_{j\ne i} u_j}d\bar{t}^2+ e^{2\alpha}(d\rho^2+dz^2)+ \sum_{j\ne i} e^{u_j} \left(d\phi^j\right)^2,
\end{split}
\end{align}
where we have used $u_i=2\log\rho-\sum_{j\ne i} u_j+c$ from \eqref{constant-lapse}, and have rescaled the time coordinate to $\bar{t}=e^{c/2}t$. Therefore, from the soliton time slice $(M^{n+2},g)$ we have obtained a vacuum spacetime $(\mathbf{M}^{n+2},\mathbf{g})$ which is static, $(n-1)$-axisymmetric, and with Kasner asymptotics. Furthermore, $\alpha$ has already been shown (Section \ref{sec4}) to be $L$-periodic in the $z$-direction. Thus, these solutions are space-periodic. Moreover, since $e^{u_i}\to 0$ on $\Gamma_{nl+i}$, these intervals on the $z$-axis become horizon rods for the Wick rotated metric. The horizon cross-sectional topologies are $S^3 \times T^{n-3}$ or $S^1 \times S^2 \times T^{n-3}$ for $n\geq 3$ depending on whether the neighboring rod structures are distinct or the same, respectively, and for $n=2$ the topology is $S^2$. Lastly, the conical singularities remain balanced on the axes as the $u_j$ and $\alpha$ have not changed. This completes the proof of Theorem \ref{wick}.  

Consider now a static solutions from either Theorem \ref{thm1} or \ref{existence}. The lapse or static potential is given by
\begin{equation}
\rho e^{-\tfrac{1}{2}\sum_{i=1}^{n}u_i}=e^{-c/2},
\end{equation}
where we have used \eqref{constant-lapse}. Since this is constant, it follows that the time slice $(M^{n+2},g)$ is Ricci flat. The remaining properties stating that this Riemannian manifold is complete, simply connected, admits a cohomogeneity-two torus action, and has Kasner asymptotics have been established above. Moreover, in the space-periodic case we have shown in Section \ref{5.2}
that the codimension-two Betti number is infinite. This completes the proof of Theorem \ref{lapse}. 
Lastly, we will establish Theorem \ref{holonomy}.

\subsection{Holonomy}
To show that the Ricci flat Riemannian manifolds $(M^{n+2},g)$ of Theorem \ref{lapse} have generic holonomy, we will make use the Ambrose-Singer Theorem~\cite{Besse}, which may be interpreted as stating that the Lie algebra $\mathfrak{hol}_p$ of the holonomy group  at $p\in M^{n+2}$ is generated by the curvature endomorphisms $R(X,Y)\colon T_q M^{n+2} \to T_q M^{n+2}$, where $X$ and $Y$ run through $T_q M^{n+2}$ and $q$ runs through $M^{n+2}$. Thus, in order to establish that the holonomy is generic, it suffices to choose a suitable $p$ and show that the curvature endomorphisms at $p$ generate $\mathfrak{so}(n+2)$. Note that the holonomy algrbra determines the holonomy group, as the manifolds discussed here are simply connected. 

We begin by computing the curvature. A straightforward although tedious computation shows that the non-zero components of the Riemann curvature tensor are given by
\begin{align}
\begin{split}
R_{\rho z \rho z} = -e^{2\alpha} \Delta_2 \alpha, \quad\quad\quad R_{\phi^i\phi^j\phi^i\phi^j} = \mathcal{G}(u_i,u_j),\quad\quad 1\leq i < j \leq n,\\
R_{\rho\phi^i\rho\phi^i} = \mathcal{F}_1 (u_i), \quad\quad
R_{\rho\phi^i z\phi^i} = \mathcal{F}_2 (u_i), \quad\quad
R_{z\phi^i z\phi^i} = \mathcal{F}_3 (u_i), \quad\quad i=1,\ldots, n,
\end{split}
\end{align}
where $\Delta_2=\partial_{\rho}^2 +\partial_z^2$ and
\begin{align}
\begin{split}
\mathcal{F}_1 (u_i) =& \frac{e^{u_i}}4 \left( 
2u_{i,\rho} \alpha_{\rho}
- 2u_{i,z} \alpha_{z}
-2 u_{i,\rho\rho}
-u_{i,\rho}^2
\right), 
\\
\mathcal{F}_2 (u_i) =& \frac{e^{u_i}}{4} \left(
2u_{i,\rho}\alpha_{z}
+ 2u_{i,z}\alpha_{\rho}
- 2 u_{i,\rho z}
- u_{i,\rho}u_{i,z}
\right), \quad\quad i=1,\ldots,n, \\
\mathcal{F}_3 (u_i) =& \frac{e^{u_i}}4 \left( 
2u_{i,z} \alpha_{z}
-2u_{i,\rho} \alpha_{\rho}
-2 u_{i,zz}
-u_{i,z}^2
\right), 
\\
\mathcal{G}(u_i,u_j)=& -\frac{e^{u_i+u_j-2\alpha}}{4} \left(
u_{i,\rho} u_{j,\rho}
+ u_{i,z}u_{j,z}
\right),\quad\quad 1\leq i<j\leq n.
\end{split}
\end{align}
It follows that the curvature endomorphisms may be represented in Weyl-Papapetrou coordinates as the matrices
\begin{equation}
    R(\partial_{\rho},\partial_z) = \begin{pmatrix}
    \begin{matrix}
    \phantom{ii}\> 0 & -e^{2\alpha}\Delta_2\alpha\\
    e^{2\alpha}\Delta_2\alpha & 0\\
    \end{matrix} 
    & \>
    \begin{matrix}
    0 & \dots & 0\\
    0 & \dots & 0
    \end{matrix} \\
    \begin{matrix}
    0 & & \> & 0 \\
    \vdots & & &  \vdots\\
    0 & & & 0
    \end{matrix} 
    &
    \begin{matrix}
    \boldsymbol{0}
    \end{matrix}
    \end{pmatrix},\quad\quad\quad\quad
   R(\partial_{\phi^i},\partial_{\phi^j}) = \begin{pmatrix}
    \begin{matrix}
    0 & 0\\
    0 & 0
    \end{matrix} 
    &
    \begin{matrix}
    0 & \dots & 0 \\
    0 & \dots  & 0
    \end{matrix}\\
    \begin{matrix}
    0 & 0 \\
    \vdots & \vdots\\
    0 & 0
    \end{matrix} 
    &
    \begin{matrix}
    \boldsymbol{G}(i,j)
    \end{matrix}
    \end{pmatrix}, 
\end{equation}
\begin{equation}
    R(\partial_{\rho},\partial_{\phi^i}) = \begin{pmatrix}
    \begin{matrix}
    0 & & & 0 \\
    0 & & & 0 \\
    \end{matrix} &
    \begin{matrix}
    \dots & \mathcal{F}_1 (u_i) & \dots & 0\\
    \dots & \mathcal{F}_2 (u_i) & \dots & 0\\
    \end{matrix} \\
    \begin{matrix}
    \vdots & \vdots \\
    -\mathcal{F}_1(u_i) & -\mathcal{F}_2(u_i) \\
    \vdots & \vdots \\
    0 & 0
    \end{matrix} &
    \begin{matrix}
    \boldsymbol{0}
    \end{matrix}
    \end{pmatrix}, 
\end{equation}
\begin{equation}
    R(\partial_z,\partial_{\phi^i}) = \begin{pmatrix}
    \begin{matrix}
    0 & & & 0 \\
    0 & & & 0 \\
    \end{matrix} &
    \begin{matrix}
    \dots & \mathcal{F}_2 (u_i) & \dots & 0\\
    \dots & \mathcal{F}_3 (u_i) & \dots & 0\\
    \end{matrix} \\
    \begin{matrix}
    \vdots & \vdots \\
    -\mathcal{F}_2(u_i) & -\mathcal{F}_3(u_i) \\
    \vdots & \vdots \\
    0 & 0
    \end{matrix} &
    \begin{matrix}
    \boldsymbol{0}
    \end{matrix}
    \end{pmatrix},
\end{equation}
where the $\mathcal{F}_k (u_i)$ (respectively $-\mathcal{F}_k (u_i)$) entries of $R(\partial_{\rho},\partial_{\phi^i})$ and $R(\partial_z ,\partial_{\phi^i})$ appear in the $(i+2)$th column (respectively $(i+2)$th row), and the only non-zero entries of the $n\times n$ matrix $\boldsymbol{G}(i,j)$ are $\mathcal{G}(u_i,u_j)$ and $-\mathcal{G}(u_i,u_j)$ which occur in its $ij$ and $ji$ entries respectively, $1\leq i<j\leq n$. Therefore, if we find a point $p\in M^{n+2}$ at which
\begin{equation}
\Delta_2\alpha\ne0, \quad\quad \mathcal{F}_1(u_i) \mathcal{F}_3(u_i) - \mathcal{F}_2(u_i)^2\ne0, \quad i=1,\dots, n, \quad\quad \mathcal{G}(u_i,u_j)\ne0,\quad 1\leq i<j\leq n,
\end{equation}
then these $(n+2)(n+1)/2$ matrices span $\mathfrak{so}(n+2)$ at that point, and the holonomy group  will be generic. 

An appropriate choice of the point $p$ is sufficiently far out in the asymptotic end. To see this, recall the asymptotics for $u_i$ and use \eqref{alpha} to obtain
\begin{equation}
u_i \thicksim \frac2n\log\rho, \quad\quad\quad \alpha \thicksim \frac{1-n}{2n}\log\rho.
\end{equation}
Furthermore, expanding in cylindrical harmonics as in \eqref{jiqughh} produces
\begin{equation}
u_{i,\rho}\thicksim\frac2{n\rho}, \quad\quad\quad \alpha_{\rho}\thicksim \frac{1-n}{2n\rho},\quad\quad\quad
u_{i,z} \text{ and } \alpha_{z} 
\text{ decay exponentially as $\rho\to\infty$}.
\end{equation}
Hence
\begin{align} \label{determinant}
\begin{split}
\mathcal{F}_1(u_i) \mathcal{F}_3(u_i) - \mathcal{F}_2(u_i)^2=&
- \frac{e^{2u_i}}{16} 
\left[
\left( 
2u_{i,\rho}\alpha_{z}
+ 2u_{i,z}\alpha_{\rho}
-2u_{i,\rho z} 
- u_{i,\rho}u_{i,z}
\right)^2 
\right.\\
-&\left. 
\left(
2u_{i,z} \alpha_{z} 
- 2u_{i,\rho} \alpha_{\rho}
-2u_{i,zz} 
-u_{i,z}^2
\right)
\left(
2u_{i,\rho} \alpha_{\rho}
-2u_{i,z} \alpha_{z} 
-2u_{i,\rho\rho} 
-u_{i,\rho}^2
\right)
\right]\\
\thicksim&
-\frac{e^{2u_i}}{8} u_{i,\rho} \alpha_{\rho}
\left(
2u_{i,\rho} \alpha_{\rho}
-2u_{i,\rho\rho} 
-u_{i,\rho}^2
\right)\\
\thicksim& \frac{(n-1)^2}{4n^4}\rho^{\frac{4}{n}-4}.
\end{split}
\end{align}
We conclude that for $\rho$ large enough, \eqref{determinant} is non-zero as $n\geq 2$. Next, observe that
\begin{equation}
\Delta_2\alpha = -\frac18 \sum_{i=1}^n|\nabla u_i|^2 +\frac{1}{2}\rho^{-2} \thicksim \frac{n-1}{2n}\rho^{-2}, \quad\quad\quad\quad
\mathcal{G}(u_i,u_j) \thicksim -\frac{\rho^{\frac{3}{n}-1}}{n^2},
\end{equation}
where the first equation is obtained from~\eqref{alpha} together with the harmonicity of $u_i$. Once again, for $\rho$ sufficiently large these expressions are non-zero. This completes the proof of Theorem~\ref{holonomy}.\bigskip\medskip


\noindent\textbf{Acknowledgements.}
The authors would like to thank Kotaro Kawai and Jordan Rainone for helpful discussions.


\end{document}